\documentclass[fontsize=10pt,twoside ]{scrreprt}

\usepackage[square,numbers,nonamebreak]{natbib}
\usepackage{setspace}

\usepackage[english]{babel}

\usepackage{chngcntr}

\usepackage[hang]{footmisc}
 at 11truept
\font\ssc=pplrc9d at 11 truept
\newcommand\qedbox{$\rlap{$\sqcap$}\sqcup$}

\usepackage[hmarginratio=1:1, bottom=2.5cm, top=2.5cm, inner=2cm, outer=2cm, paperwidth=17cm, paperheight=24cm]{geometry}

\usepackage[automark]{scrlayer-scrpage}
\cohead{\textnormal{Fundamental Functor on Hypergroups}}
\lehead{\pagemark}
\rohead{\pagemark}

\let\ceheadL\cehead
\renewcommand\cehead[1]{
\ceheadL{\textnormal{#1}}
}

\cehead{Behnam Afshar --- Reza Ameri}
\lefoot{}
\rofoot{}

\usepackage{graphicx}
\usepackage[usenames,dvipsnames,svgnames,table]{xcolor}
\definecolor{Maroon}{cmyk}{0, 0.87, 0.68, 0.32}
\definecolor{RoyalBlue2}{cmyk}{80,100,0,0.1}

\newcommand\auths[1]{\large \textsc{\textcolor{Maroon}{#1}}\setstretch{1.2}}
\newcommand\titl[1]{\center \linespread{1.1}\color{RoyalBlue2}\Large\textbf{ #1}\color{black}\bigskip} 
\renewcommand\abstract[1]{
\begin{center}
{\textbf{Abstract}}
\end{center}
{
\linespread{1.1}\fontsize{9pt}{-10pt}\selectfont #1}}
\usepackage{tikz}
\usepackage{tikz-cd} 

\usepackage[utf8]{inputenc}
\usepackage[sc,osf]{mathpazo}
\usepackage[euler-digits]{eulervm}
\usepackage[T1]{fontenc}
\usepackage{mathtools}

\DeclareSymbolFont{operators}{\encodingdefault}{ppl}{m}{n}
\DeclareMathAlphabet{\mathbf}{\encodingdefault}{ppl}{bx}{n}
\DeclareMathAlphabet{\mathit}{\encodingdefault}{ppl}{m}{it}

\usepackage{amsfonts}
\usepackage{amsmath}
\usepackage{ragged2e}
\usepackage{titlesec}

\renewcommand{\thesection}{\arabic{section}}
 
\titleformat{\section}{\medskip\bigskip\normalfont\Large\bf}{\thesection}{0.5em}{}
\titleformat{\subsection}{\smallskip\bigskip\normalfont\large\bf}{\thesubsection}{0.5em}{}

\usepackage{amsthm}
\newtheoremstyle{dotless}{}{}{\itshape}{}{\bfseries}{}{1em}{}

\theoremstyle{dotless}
\newtheorem{theo}{Theorem}
\newtheorem{prop}[theo]{Definition}
\newtheorem{lem}[theo]{Lemma}
\newtheorem{cor}[theo]{Corollary}
\newtheorem{rem}[theo]{Remark}
\newtheorem{ex}[theo]{Example}
\renewenvironment{proof}{\smallbreak\noindent {\sc Proof \;---\;}}{\hfill\qedbox}

\numberwithin{theo}{section}
\numberwithin{equation}{section}

\DeclareOldFontCommand{\rm}{\normalfont\rmfamily}{\mathrm}
\DeclareOldFontCommand{\sf}{\normalfont\sffamily}{\mathsf}
\DeclareOldFontCommand{\tt}{\normalfont\ttfamily}{\mathtt}
\DeclareOldFontCommand{\bf}{\normalfont\bfseries}{\mathbf}
\DeclareOldFontCommand{\it}{\normalfont\itshape}{\mathit}
\DeclareOldFontCommand{\sl}{\normalfont\slshape}{\@nomath\sl}
\DeclareOldFontCommand{\sc}{\normalfont\scshape}{\@nomath\sc}

\usepackage[hidelinks]{hyperref}
\begin{document}

\setheadsepline{1pt}[\color{black}]

\hskip 1 cm

\titl{Fundamental Functor on Hypergroups
}

\auths{Behnam Afshar --- Reza Ameri}

\center{{\fontsize{9pt}{-10pt}{\it\selectfont Department of Mathematics, Statistics and Computer Science,\\ College of Science, University of Tehran, Tehran, Iran}}}

\thispagestyle{empty}
\justify\noindent
\setstretch{0.3}
\abstract{For a hypergroup $(H,\circ)$ we consider $\gamma^{\ast}$, as the smallest  equivalence relation on $H$ such that the quotion $(H/\gamma^{\ast},\tiny{\otimes})$ is an abelian group. We study some more properties of $\gamma^{\ast}$.  Initially, it is investigated which subhypergroup the congruence relation modulo is strongly regular on, and its quotient results in an abelian group? This is directly related to the fundamental relation $\gamma^{\ast}$, since such subhypergroups must contain $S_{\gamma}$. Then, we examine the functor $\gamma^{\ast}$ from a categorical perspective and investigate properties such as continuity and cocontinuity concerning it using the decomposition $\gamma=\delta\tiny{\ast}\beta$. For this purpose, we define the reduced words on strongly regular hypergroups. This has a direct application in studying how the functor $\gamma^{\ast}$ affects on the stalks of the sheaves of hypergroups.}

\setstretch{2.1}
\noindent
{\fontsize{10pt}{-10pt}\selectfont {\it Mathematics Subject Classification (2010)}: 20N20, 18D35, 35A27}\\[-0.8cm]

\noindent 
\fontsize{10pt}{-10pt}\selectfont  {\it Keywords}: derived  subhypergroup; strongly regular relation; reduced word\\[-0.8cm]

\setstretch{1.1}
\fontsize{11pt}{12pt}\selectfont

\section{Introduction}

On of the most important topics in hyperstructure theory is the relationship between hyperstructures and their classical state, which is done with fundamental relations. The fundamental relation $\beta^{\ast}$, which is the transitive closure of $\beta$ was introduced by M. Koskas \cite{Koskas} and studied by P. Corsini \cite{Corsini}, D. Freni \cite{Freni1,Freni2} and T. Vogiouklis \cite{Vougiouklis}. Then Freni introduced the fundamental relation $\gamma^{\ast}$ which is the transitive closure of $\gamma$ and is the smallest strongly regular relation on a given hypergroup $H$ such that $H/\gamma^{\ast}$ is an abelian group. He also showed that the relations $\beta$ and $\gamma$ are transitive in hypergroups. J. Jantosciak studied the relationship between regular relations and subhypergroups of a transposition hypergroup. He proved that there is a one-to-one correspondence between regular relations and reflexive closed subhypergroups of transposition hypergroup \cite{Jantosciak}.
In \cite{AmeriAfshar1} it is proved that $\gamma^{\ast}=\delta\ast\beta^{\ast}$, where $\delta$ is a congruence relation with respect to the commutator subgroup and, in \cite{AmeriAfshar2} it has been shown that strongly regular relations on hypergroups can be considered as congruence relations. In other words, it is proved in \cite{AmeriAfshar2} that there exists a correspondence relation between the strongly regular relations and  normal subhypergroups containing $S_{\beta}$ the heart of $\beta$.\\
In addition, a covariant functor from the category of semihypergroups to the category of semigroups (resp. abelian semigroups) is obtained from the fundamental relation $\beta^{\ast}$, (resp. $\gamma^{\ast}$) \cite{DavvazFotea}. Another important concept is complete parts, which was introduced by M. Koskas \cite{Koskas} and, studied by P. Corsini \cite{Corsini} and Y. Sureau \cite{Sureau} in the hypergroup theory. In \cite{AmeriAfshar2} it is proved that any normal subhypergroup containing $S_{\beta}$ is complete part.
In recent years, the study of hyperstructures from the perspective of category theory has gained attention. Through this, important categories such as hypergroups, hyperrings, hypermodules, ect., have been investigated, and their relationships with classical categories have been studied \cite{Ameri,AmeriShojaei,NorouziAmeri}.\\
In Section \ref{sec1}, we examine some properties of the relation $\gamma^{\ast}$, which is actually a generalization for the results of \cite{AmeriAfshar2} and in Section \ref{sec2}, we will study properties such as continuity and co-continuity for $\gamma^{\ast}$ as a functor from the category of strongly regular hypergroups to the category of abelian groups. The main use of this research is that it helps us to understand how $\gamma^{\ast}$ acts on stalks of the sheaves of hypergroups.\\

\section{Preliminaries}
Let $(H,\circ)$ be a non-empty set and $\circ:H\tiny{\times}H\rightarrow \mathcal{P}^*(H)$ be a hyperoperation where $\mathcal{P}^*(H)$ is the family of non-empty subsets of $H$. The couble $(H,\circ)$ is called a \textit{hypergroupoid}. For any two non-empty subsets $A$ and $B$ of $H$ we define $A\tiny{\circ} B=\bigcup_{(a,b)\in A\tiny{\times}B}{a\tiny{\circ} b}$.\\
A hypergroupoid $(H,\circ)$ is called a \textit{semi-hypergroup} if for all $a,b,c$ of $H$ we have $(a\tiny{\circ}b)\tiny{\circ}c=a\tiny{\circ}(b\tiny{\circ}c)$. A hypergroupoid $(H,\circ)$ is called a \textit{quasi-hypergroup} if for all $a$ of $H$ we have  $a\tiny{\circ}H=H\tiny{\circ}a=H$.\\
A  hypergroupoid $(H,\circ)$ which is both a semi-hypergroup and a quasi-hypergroup, is called a \textit{hypergroup}.
Let $(H,\circ)$ be a semi-hypergroup and $\mathcal{R}$ be an equivalence relation on $H$. If $A$ and $B$ are non-empty subsets of $H$ then $A\bar{\mathcal{R}}B$ means that for each $a\in A$, there exist $b\in B$ such that $a\mathcal{R}b$ and for each $b'\in B$ there exist $a'\in A$ such that $a'\mathcal{R}b'$. $A \bar{\bar{\mathcal{R}}}B$ means that for each $a\in A$ and $b\in B$ we have $a\mathcal{R}b$ \cite{Survey2}.\\

{\prop\emph{\cite{Survey2}}
	Let $(H,\circ)$ be a semi-hypergroup. The equivalence relation $\mathcal{R}$ is called:\\
1- Regular on the right (on the left), if for all $x\in H$, from $a\mathcal{R}b$, it follows that $(a\tiny{\circ}x)\bar{\mathcal{R}}(b\tiny{\circ}x)$ $((x\tiny{\circ} a)\bar{\mathcal{R}}(x\tiny{\circ}b) \ respectively)$;\\
2- Strongly regular on the right (on the left), if for all $x\in H$, from $a\mathcal{R}b$, it follows that $(a\tiny{\circ} x)\bar{\bar{\mathcal{R}}}(b\tiny{\circ}x)$ $((x\tiny{\circ} a)\bar{\bar{\mathcal{R}}}(x\tiny{\circ}b) \ respectively)$;\\
3- Regular (strongly regular), if it is regular (strongly regular) on the right and on the left.}\\

If $\mathcal{R}$ is a strongly regular relation on the hypergroup $H$, throughout this paper, we denote the group $H/\mathcal{R}$ by the symbol $\mathcal{R}(H)$ and  for any normal subhypergroup $K$ of $H$; $H/K=\{h\circ K ; h\in H\}$.

{\theo\emph{\cite{Survey2}}
	Let $(H,\circ)$ be a semi-hypergroup and $\mathcal{R}$ be an equivalence relation on $H$;\\
	If $\mathcal{R}$ is regular, then $H/\mathcal{R}$ is a semi-hypergroup with respect to the hyperoperation $\mathcal{R}(x)\tiny{\otimes}\mathcal{R}(y)=\{{\mathcal{R}(z);  z\in x\tiny{\circ}y}\}$;\\
	If the above hyperoperation is well defined on $H/\mathcal{R}$, then $\mathcal{R}$ is regular.}

{\cor\emph{\cite{Survey2}}
	If $(H,\circ)$ is a hypergroup and $\mathcal{R}$ be an equivalence relation on $H$, then $\mathcal{R}$ is regular (strongly regular) if and only if $(\mathcal{R}(H),\otimes)$ is a hypergroup (group).}

{\theo\emph{\cite{Survey2}}
	Let $(H,\circ)$ be a semi-hypergroup and $\mathcal{R}$ be an equivalence relation on $H$.\\
	If $\mathcal{R}$ is strongly regular, then $H/\mathcal{R}$ is a semi-group with respect to the operation $\mathcal{R}(x)\tiny{\otimes}\mathcal{R}(y)=\mathcal{R}(z)$, for all $z\in x\tiny{\circ}y$;\\
	If the above operation is well defined on $H/\mathcal{R}$, then $\mathcal{R}$ is strongly regular.}\\

Let $(H,\circ)$ be a semihypergroup and $C$ be a non-empty
subset of $H$. We say that $C$ is a \textit{complete part} of $H$ if for any non-zero natural number $n$ and for all $a_{1}, . . . , a_{n}$ of $H$, the following implication holds:
\begin{center}
	$C\cap\prod_{i=1}^{n}a_{i}\neq\emptyset\Rightarrow\prod_{i=1}^{n}a_{i}\subseteq C.$
\end{center}

{\prop\emph{\cite{Survey1}}
	For all $n>1$, we define the relations $\beta_{n}$ and $\gamma_n$ on a semi-hypergroup $H$, as follows:
	\begin{center}
		$a\beta_{n}b\Longleftrightarrow \exists(x_{1},x_{2},...,x_{n})\in H^{n} : \{a,b\}\subseteq\prod_{i=1}^{n}x_{i}$
	\end{center}
	\begin{center}
		$ a\gamma_{n}b\Longleftrightarrow \exists(x_{1},x_{2},...,x_{n})\in H^{n} , \sigma\in \mathbb{S}_{n}: \ a\in\prod_{i=1}^{n}x_{i}  \ , \  b\in\prod_{i=1}^{n}y_{\sigma(i)}$
	\end{center}
	and $\beta=\bigcup_{n\geq1}\beta_{n}$ and $ \gamma=\bigcup_{n\geq1}\gamma_{n}$ where $\beta_{1}=\gamma_{1}=\{(x,x);x\in H\}$.
	Let $\beta^{*}$be the transitive closure of $\beta$ and $\gamma^{*}$ be the transitive closure of $\gamma$.}\\

Both relations $\beta^{*}$ and $\gamma^{*}$ are strongly regular and transitive in hypergroups, i.e. $\gamma=\gamma^{*}$ and $\beta=\beta^{*}$. If $H$ is a hypergroup then $\beta^{*}(H)$ is called the \textit{fundamental group} and $\beta^{*}$ is the smallest strongly regular relation on $H$. 
The relation $\gamma^{*}$ is the smallest strongly regular relation on a semi-hypergroup $H$ (resp. hypergroup $H$) such that the quotient $\gamma^{*}(H)$ is commutative semi-group (resp. commutative group), \cite{Corsini}.

{\prop\emph{\cite{Freni1}}
	Let $D_{1}$, $D_{2}$ and $D$ denote the sets
	\begin{center}
		$D_{1}=\underset{(x,y)\in H^{2}}{\bigcup}xy/yx$,\ \ \  $D_{2}=\underset{(x,y)\in H^{2}}{\bigcup}xy\backslash yx$, \ \ \ $D=D_{1}\cup D_{2}.$
	\end{center}
	Then the derived subhypergroup $H'$ defined as the intersection of all subhypergroups that are complete pare and contain $D$.}\\

A subhypergroup $K$ of $H$ is called \textit{conjugable on the right} (\textit{on the left}) if for all $k_{1}$, $k_{2}$ of $K$ and $x$ of $H$,
from $k_{1}\in k_{2}\circ x$ ($k_{1}\in x\circ k_{2}$, respectively), it follows that $x\in K$
and there exists $x'\in H$ such that $x'\circ x\subseteq K$. We say that $K$ is \textit{conjugable} if it is conjugable on the left and on the right \cite{DavvazFotea}.\\

Let $SR(H)$ to be the set of all strongly regular relations on hypergroup $H$ and for every $\rho\in SR(H)$, $S_{\rho}:=\{x\in H \ ;\ \rho(x)=e_{\rho(H)}\}$, and $S(H)=\{S_{\rho} \ ;\ \rho\in SR(H)\}$.
Clearly, $S_{\rho}$ is a complete part and conjugable subhypergroup  \cite{AmeriAfshar2}. We will also display the set of all subhypergroups containing $S_{\beta}$ and the set of all normal subhypergroups containing $S_{\beta}$ with symbols $(S_{\beta})$ and $N(S_{\beta})$, respectively. That, $S_{\beta}$ is actually the heart of $H$ and the following theorem is proved about it.

{\theo\emph{\cite{DavvazFotea}}\label{q4}
	The heart of a hypergroup $H$ is the intersection of all subhypergroups of $H$ which are complete parts.}

{\theo\emph{\cite{AmeriAfshar2}}
	If $H$ is a regular hypergroup then:\\
	(i) $K\in N(S_{\beta})$ if and only if $H/K$ is group;\\
	(ii) $K\in (S_{\gamma})$ if and only if $H/K$ is abehian group.}\\

For regular hypergroup $H$, we have:
\begin{center}
	$(S_{\gamma})=N(S_{\gamma})$ \ \ ,\ \  $ N(S_{\beta})=S(H)\cong SR(H).$ 
\end{center} 
If $K\in(S_{\gamma})$, then $K$ is complete part, conjugable and normal subhypergroup of $H$, \cite{AmeriAfshar2}.\\
Let $\rho$ be a strongly regular relation on hypergroup $H$, and $\sigma$ be a congruence relation on group $\rho(H)$, then consider:
\begin{center}
	$(a,b)\in\sigma{\tiny{\ast}}\rho\Longleftrightarrow(\rho(a),\rho(b))\in\sigma.$
\end{center} 
If $\delta$ is congruence relation on $G$ containing module the commutator subgroup of $G$, then $\delta{\tiny{\ast}}\rho$ is the smallest strongly regular relation containing $\rho$ on $H$, such that $H/\delta{\tiny{\ast}}\rho$ is commutative group. In particular, $\gamma=\delta{\tiny{\ast}}\beta$ \cite{AmeriAfshar1}.

{\prop\emph{\cite{Awodey}}
	A functor $F:\textbf{C}\rightarrow\textbf{D}$ is said to preserve limits of type $\textbf{J}$ if,
	whenever $p_{j}:L\rightarrow D_{j}$ is a limit for a diagram $D: \textbf{J}\rightarrow \textbf{C}$; the cone $Fpj : FL\rightarrow FDj$ is then a limit for the diagram $FD:\textbf{J}\rightarrow \textbf{D}$. Briefly
	\begin{center}
		$F(\underset{\textbf{J}}{lim}D_{j})\cong\underset{\textbf{J}}{lim}F(D_{j})$.
	\end{center}
	A functor that preserves all limits is said to be continuous. Similarly, the cocontinuous functor is also defined through the concept of colimit.}

\section{Canonical quotient hypergroups and $\gamma^{\ast}$ relation} \label{sec1}

Let $H'$ be derived subhypergroup of $H$ and $S_{\rho}\in S(H)$. Then $S_{\gamma}=H'$, \cite{Freni1} and, $S_{\rho}$ is a conjugable, normal and complete part subhypergroup of $H$, containing $S_{\beta}$, \cite{AmeriAfshar2}.

{\lem\label{1111}
	For every $\rho\in SR(H)$ and $x\in H$; $\rho(x)=x\circ S_{\rho}$.}
{\begin{proof}	
		If $y\in x\circ S_{\rho}$, then $\rho(x)=\rho(y)$ and $x\circ S_{\rho}\subseteq\rho(x)$. Let $y\in \rho(x)$. Since $H=x\circ H$, then there exist $h\in H$ such that $y\in x\circ h$. So $\rho(y)=\rho(x)\tiny{\otimes}\rho(h)$ and $\rho(x)=\rho(y)$. Therefore $h\in S_{\rho}$ and $y\in x\circ S_{\rho}$.		
	\end{proof}}


{\theo   \label{2222}
	If $K$ is a closed subhypergroup of $H$, then $H/K$ is group if and only if  $K\in N(S_{\beta})$.}
	{\begin{proof}
		If $H/K$ is a group and $(x_{1},...,x_{n})\in H^{n}$, then:
		\begin{center}
			$\prod_{i=1}^{n}x_{i}\circ K= (\prod_{i=1}^{n}x_{i})\circ K=\{t\circ K; t\in \prod_{i=1}^{n}x_{i}\}= r\circ K$
		\end{center}
		where $r\in \prod_{i=1}^{n}x_{i}$. If $\prod_{i=1}^{n}x_{i}\bigcap K\neq\emptyset$, then there exist $r\in \prod_{i=1}^{n}x_{i}\bigcap K$ such that $(\prod_{i=1}^{n}x_{i})\circ K= r\circ K= K$. Since $K$ is closed, then $\prod_{i=1}^{n}x_{i}\subseteq K$.  Thus $K$ is complete part subhypergroup of $H$ and since $S_{\beta}$ is the smallest complete part subhypergroup of $H$ (see Theorem \ref{q4}), then $S_{\beta}\subseteq K$. Therefore $K\in N(S_{\beta})$.\\
		Let $K\in N(S_{\beta})$. Since $K\lhd H$, then $H/K$ is hypergroup. From $S_{\beta}\subseteq K$, it follows that the hyperoperation $\cdot:H/K\times H/K\rightarrow \mathcal{P}^{\ast}(H/K)$ where $(x\circ K)\cdot(y\circ K)=(x\circ y)\circ K$ is an operation. Because if $(x\circ y)\bigcap K\neq\emptyset$, then for every $t\in (x\circ y)\bigcap K$ we have $(x\circ y)\circ K= (x\circ y)\circ S_{\beta}\circ K= t\circ S_{\beta}\circ K= t\circ K= K$. So $x\circ y \subseteq K$. Also if $(x\circ y)\bigcap K=\emptyset$, then for every $t\in x\circ y$; $(x\circ y)\circ K= (x\circ y)\circ S_{\beta}\circ K= t\circ S_{\beta}\circ K= t\circ K$.
	\end{proof}}

\bigskip

Let $K$ be a closed subhypergroup of $H$ and $S_{\gamma}\subseteq K$. Then for every $x\in H$, there exists $x'\in H$ such that $\gamma(x)\tiny{\otimes}\gamma(x')=\gamma(x\circ x')=\gamma(x'\circ x)=e_{\gamma(H)}$. For every $y\in K$;
\begin{center}
	$(x\circ y \circ x')\circ S_{\gamma}=(x'\circ y \circ x)\circ S_{\gamma}=(y\circ x \circ x')\circ S_{\gamma}=y\circ S_{\gamma}$.
\end{center}
So $(x\circ y \circ x')\circ S_{\gamma}\circ K=(x'\circ y \circ x)\circ S_{\gamma}\circ K=y\circ S_{\gamma}\circ K$. Then $(x\circ y \circ x')\circ K= K$ and since $K$ is closed, then $x\circ K\circ x'\subseteq K$ for every $x\in H$, and similarly $x'\circ K\circ x\subseteq K$. Therefore $x\circ K\subseteq K\circ x$ and $K\circ x\subseteq x\circ K$. Hence $K\circ x=x\circ K$.

{\cor
	If $K$ is a closed subhypergroup of hypergroup $H$, then $H/K$ is abelian group if and only if  $K\in (S_{\gamma})$.}
	{\begin{proof}
		Since $S_{\gamma}$ is the smallest subhypergroup of $H$ such that $H/S_{\gamma}$ is abelian group, then $S_{\gamma}\subseteq K$. Let $K\in (S_{\gamma})$, since $K$ is closed, then $K\in N(S_{\beta})$ and by Theorem \ref{2222} $H/K$ is group and if $x\circ K, y\circ k \in H/K$, then:
		\begin{align*}
			(x\circ K)(y\circ K)&=(x\circ y)\circ K= (x\circ y)\circ S_{\gamma}\circ K= (y\circ x)\circ S_{\gamma}\circ K\\
			&=(y\circ x)\circ K= (y\circ K)(x\circ K).	\end{align*}\end{proof}}

{\theo
	Every regular relation containing $\beta$ is strongly regular.}	
{	\begin{proof}	
		Let $\rho$ be a regular relation on $H$ and $\beta\subseteq\rho$. For every $x,y,z\in H$, if $(x,y)\in\rho$, then $x\circ z\ \bar{\rho}\ y\circ z$ and since $\beta\subseteq\rho$, then $(x_{1},x_{2}),(y_{1},y_{2})\in\beta\subseteq\rho$ for every $x_{1},x_{2}\in x\circ z$ and $y_{1},y_{2}\in y\circ z$. So $x\circ z\ \bar{\bar{\rho}} \ y\circ z$.
	\end{proof}}

\bigskip

A subhypergroup of a canonical hypergroup, is called a \textit{canonical subhypergroup} if it is a canonical hypergroup itself and with the same identity element. Also a subhypergroup $K$ of canonical hypergroup $H$ is canonical subhypergroup if and only if $0_{H}\in K$, \cite{Massouros}. Therefore $S_{\rho}$ is canonical subhypergroup of $H$, for every $\rho\in SR(H)$.
If $K$ is a canonical subhypergroup of canonical hypergroup $H$, then $H/K$ is canonical hypergroup \cite{Massouros}.

{\cor	If $N$ is a canonical subhypergroup of canonical hypergroup $H$, then $\overset{N}{\equiv}$ is regular.}
	{\begin{proof}
		Let $x\circ N=y\circ N$	and $z\in H$. So $(x\circ z)\circ N=(y\circ z)\circ N$ and if $a\in x\circ z$, then $a\in(x\circ z)\circ N=(y\circ z)\circ N$. So there are $b\in y\circ z$ and $n\in N$ such that $a\in b\circ n$, thus $a\circ N\subseteq b\circ N$. Also $b\in a\circ n^{-1}$, so $b\circ N\subseteq a\circ N$ and therefore $a\circ N=b\circ N$.
	\end{proof}}

\bigskip

It is easy to prove that, if $\rho\in SR(H)$, $\sigma=\ \overset{K}{\equiv}$ and $\tau\in SR(H/K)$ where $K$ is a normal subhypergroup of canonical hypergroup $H$, then $\rho\vee\sigma, \tau\ast\sigma\in SR(H)$.

{\theo\label{6666}
	If $K$ is a canonical subhypergroup of canonical hypergroup $H$, then:
	\begin{equation}
		S(H/K)\cong N(\dfrac{S_{\beta}\tiny{\circ}K}{K}).
	\end{equation}}
{\begin{proof}Clearly $\boldsymbol\beta(h\circ K)=\beta(h)\circ K$, and by Lemma \ref{1111}, $\boldsymbol\beta(h\circ K)=h\circ S_{\beta}\tiny{\circ}K$. Then:
		\begin{align*}
			S_{\boldsymbol\beta} 
			&= \{h\circ K;\  \boldsymbol\beta(h\circ K)\ .\ \boldsymbol\beta(x\circ K)=\boldsymbol\beta(x\circ K),\  \forall x\circ K\in H/K\}\\
			&= \{h\circ K;\ h\circ S_{\beta}\tiny{\circ}K \ . \ x\circ S_{\beta}\tiny{\circ}K=x\circ S_{\beta}\tiny{\circ}K, \ \forall x\in H\}\\
			&= \{h\circ K;\  (h\circ x)\circ S_{\beta}\tiny{\circ}K=x\circ S_{\beta}\tiny{\circ}K,\ \forall x\in H\}\\
			&= \{h\circ K;\ h\in S_{\beta}\tiny{\circ}K\}\\
			&= \dfrac{S_{\beta}\tiny{\circ}K}{K}.	
		\end{align*}\end{proof}}

Just as there is a one-to-one correspondence between the set of subhypergroups of canonical hypergroup $H$ containing canonical subhypergroup $K$, and the set of canonical subhypergroups of $H/K$, there is also a one-to-one correspondence between the set of canonical subhypergroups of $H$ containing $S_{\beta}\tiny{\circ}K$ and the set of canonical subhypergroups of $H/K$ containing $S_{\boldsymbol\beta}=\dfrac{S_{\beta}\tiny{\circ}K}{K}$, where $S_{\boldsymbol\beta}=\omega_{H/K}$ and $S_{\beta}=\omega_{H}$.\\
Therefore, the set of strongly regular relations on $H$ containing $\overset{K}{\equiv}$, is in one-to-one correspondence with the set of strongly regular relations on $H/K$. Let $\sigma=\ \overset{K}{\equiv}$, then:
\begin{align*}
	SR(H/K) 
	&\cong SR(H)\vee\sigma=\{\rho\vee\sigma;\ \rho\in SR(H)\};\\
	S(H/K)
	&\cong S(H)\tiny{\circ}K=\{S_{\rho}\tiny{\circ}K;\ \rho\in SR(H)\}.	
\end{align*}

{\rem
	If $\boldsymbol\rho\in SR(H/K)$, because $\boldsymbol\rho(h\circ K)=\rho(h)\circ K= h\circ S_{\rho}\tiny{\circ}K$, then with a proof similar to Theorem \ref{6666} we will have $S_{\boldsymbol\rho}=\dfrac{S_{\rho}\tiny{\circ}K}{K}$ and therefore:
	\begin{center}
		$S(H/K)= N(\dfrac{S_{\beta}\tiny{\circ}K}{K})=\Big\{\dfrac{S_{\rho}\tiny{\circ}K}{K};\ \rho\in SR(H)\Big\}$.
	\end{center}}

{\theo
	If $K$ is a canonical subhypergroup of canonical hypergroup $H$ and  $\boldsymbol\rho\in SR(H/K)$, then:
	\begin{equation}
		\boldsymbol\rho(H/K)=\dfrac{H}{S_{\rho}\tiny{\circ}K}.
	\end{equation}}
	{\begin{proof}		
		Since 	$\boldsymbol\rho(H/K)=\dfrac{H/K}{S_{\boldsymbol\rho}}$ and $S_{\boldsymbol\rho}=\dfrac{S_{\rho}\tiny{\circ}K}{K}$, then by third isomorphism theorem $\boldsymbol\rho(H/K)=\dfrac{H}{S_{\rho}\tiny{\circ}K}$.\end{proof}}

\bigskip

Thus the following diagram is commutative:
	\begin{center}
	\begin{tikzcd}	
		H \ar[r,"\rho"]  \ar[rd,"\rho\vee\sigma"] \ar[d,"\sigma"']  &\dfrac{H}{S_{\rho}} \ar[d,"\sigma'"] \\
		\dfrac{H}{K}   \ar[r,"\boldsymbol\rho"']	& \dfrac{H}{S_{\rho}\tiny{\circ}K}
	\end{tikzcd}
\end{center}
Where $\sigma$ and $\sigma'$ are congruence relations modulo $K$ and $\dfrac{S_{\rho}\tiny{\circ}K}{S_{\rho}}$ respectively. Therefore:
\begin{equation} 
	\boldsymbol\rho\ast\sigma=\rho\vee\sigma=\sigma'\ast\rho.
\end{equation}
As a final result $S_{\boldsymbol\gamma}=\big(\dfrac{H}{K}\big)'=\dfrac{S_{\gamma}\tiny{\circ}K}{K}= \dfrac{H'\tiny{\circ}K}{K}$ and $\boldsymbol\gamma(\dfrac{H}{K})=\dfrac{\frac{H}{K}}{S_{\boldsymbol\gamma}}=\dfrac{H}{H'\tiny{\circ}K}$. In addition
\begin{equation} 
	(\sigma'\ast\gamma)(H)=\sigma'\big(\dfrac{H}{H'}\big)=\dfrac{\frac{H}{H'}}{\frac{H'\tiny{\circ}K}{H'}}=\dfrac{H}{H'\tiny{\circ}K}.
\end{equation}

{\ex
	Let $H$ be the following canonical hypergroup.	
	\begin{center}	
		\begin{tabular}{c|ccccccccc}
			$\circ$&$e$& $a$ & $b$ & $c$ & $x$ & $y$ & $z$ & $u$ & $v$\\
			\hline
			$e$ & $e$ &$a$ &$b$ &$c$ &$x$   &$y$   &$z$    &$u$   &$v$\\
			$a$ & $a$ &$e$ &$c$ &$b$ &$x$   &$y$   &$u$    &$z$   &$v$\\
			$b$ & $b$ &$c$ &$e$ &$a$ &$y$   &$x$   &$z$    &$u$   &$v$\\
			$c$ & $c$ &$b$ &$a$ &$e$ &$y$   &$x$   &$u$    &$z$   &$v$\\
			$x$ & $x$ &$x$ &$y$ &$y$ &$b,c$ &$e,a$ &$v$    &$v$   &$z,u$\\
			$y$ & $y$ &$y$ &$x$ &$x$ &$e,a$ &$b,c$ &$v$    &$v$   &$z,u$\\
			$z$ & $z$ &$u$ &$z$ &$u$ &$v$   &$v$   &$a,c$  &$e,b$ &$x,y$\\
			$u$ & $u$ &$z$ &$u$ &$z$ &$v$   &$v$   &$e,b$  &$a,c$ &$x,y$\\
			$v$ & $v$ &$v$ &$v$ &$v$ &$z,u$ &$z,u$ &$x,y$  &$x,y$ &$e,a,b,c$		
		\end{tabular}
	\end{center}
	Then $K=\{e,a\}$ is a canonical subhypergroup of $H$ and $b\circ K=c\circ K=\{b,c\}$, $x\circ K=\{x\}$, $y\circ K=\{y\}$, $z\circ K=u\circ K=\{z,u\}$ and $v\circ K=\{v\}$. Thus $H/K=\{K, b\circ K, x\circ K, y\circ K, z\circ K, v\circ K \}$ is the following canonical hypergroup.
	\begin{center}
		\begin{tabular}{c|cccccc}
			
			$\ast$&   $K$& $b\circ K$ & $x\circ K$ & $y\circ K$ & $z\circ K$ & $v\circ K$\\
			\hline
			$K$       & $K$       & $b\circ K$ & $x\circ K$ &$y\circ K$  &$z\circ K$     & $v\circ K$\\
			$b\circ K$& $b\circ K$& $K$        &$y\circ K$  &$x\circ K$  &$z\circ K$     &$v\circ K$\\
			$x\circ K$& $x\circ K$& $y\circ K$ &$b\circ K$  &$K$         &$v\circ K$   &$z\circ K$\\
			$y\circ K$& $y\circ K$& $x\circ K$ &$K$         &$b\circ K$  &$v\circ K$   &$z\circ K$\\
			$z\circ K$& $z\circ K$& $z\circ K$ &$v\circ K$  &$v\circ K$  &$K$,$b\circ K$   &$x\circ K$,$y\circ K$\\
			$v\circ K$& $v\circ K$& $v\circ K$ &$z\circ K$  &$z\circ K$  &$x\circ K$,$y\circ K$   &$K$,$b\circ K$
			
		\end{tabular}
	\end{center}
	Clearly $\beta=K^{2}\cup\{(x,y),(y,x),(z,u),(u,z),(v,v)\}\cup\Delta_{H}$ and $S_{\beta}=\{e,a,b,c\}$. Also $\beta(H)=\{\{e,a,b,c\},\{x,y\},\{z,u\},\{v\}\}\cong V_{4}$.\\
	On the other hand $\boldsymbol\beta=\{(K,b\circ K),(b\circ K,K)\}\cup\{(y\circ K,x\circ K),(x\circ K,y\circ K)\}\ \cup\Delta_{H/K}$ and $\boldsymbol\beta(H/K)=\{\{K,b\circ K\},\{x\circ K,y\circ K\},\{z\circ K\},\{v\circ K\}\}\cong V_{4}$. So $S_{\beta}\circ K=\{e,a,b,c\}$ and $S_{\boldsymbol\beta}=\dfrac{S_{\beta}\circ K}{K}=\{K,b\circ K\}$.}

\section{Categorical properties of $\gamma^{*}$}\label{sec2}

If $H$, $K$ are hypergroups, then $S_{\beta_{H\times K}}=S_{\beta_{H}}\times S_{\beta_{K}}$ and $\beta(H\times K)=\frac{H\times K}{S_{\beta_{H\times K}}}=\frac{H}{S_{\beta_{H}}}\times\frac{K}{S_{\beta_{K}}}=\beta(H)\times\beta(K)$.
On the other hand, if $G_{1}$ and $G_{2}$ are groups, then $(G_{1}\times G_{2})'=G_{1}'\times G_{2}'$. Then:
\begin{center} 
	$\gamma(H\times K)=\frac{\beta(H\times K)}{\beta(H\times K)'}=\frac{\beta(H)\times\beta(K)}{\beta(H)'\times\beta(K)'}=
	\frac{\beta(H)}{\beta(H)'}\times\frac{\beta(K)}{\beta(K)'}=\gamma(H)\times \gamma(K)$.
\end{center} 
Hence, if $\{H_{i}\}_{i\in I}$ is a family of hypergroups, then $\gamma(\prod_{i\in I}H_{i})=\prod_{i\in I}\gamma(H_{i})$. We want to check this, for coproduct and colimit as well.

{\theo The functors $\delta$ and $\beta$ are continuous.}
	{\begin{proof} Consider the diagram
		\begin{center}
			\begin{tikzcd}	
				& Group \ar[d,"\delta"] \\
				I \ar[r,"\delta F"']  \ar[ur,"F"]	& A.Group
			\end{tikzcd}
		\end{center}
		and suppose that $\underset{I}{lim} \mathcal{F}_{j}$ together with a family $\{p_{i}:\underset{I}{lim} \mathcal{F}_{j}\rightarrow F_{i} \}_{i\in I}$ of morphisms is limit of $I$-diagram $F$. Let $G=\underset{I}{lim} \mathcal{F}_{j}$, then we must show that $\delta(\underset{I}{lim} \mathcal{F}_{j})=G/G'$ together with the family $\{\delta p_{i}:G/G'\rightarrow G_{i}/G'_{i}\}_{i\in I}$ of morphisms is also a limit of $I$-diagram $\delta F$.\\
		(i) Since $G$ is a limit of $I$-diagram $G$, then for every $\phi:j\rightarrow i$ in $I$, we have $P_{i}=\mathcal{G}(\phi)\circ P_{j}$. So $\delta P_{i}=\delta \mathcal{G}(\phi)\circ\delta P_{j}$.\\
		(ii) Let $A$ be an object of category $A.Group$ and $\{q_{i}:A\rightarrow G_{i}/G'_{i}=\delta(\mathcal{G}_{i})\}_{i\in I}$ be a family of morphisms in $A.Group$ such that for every $\phi:j\rightarrow i$ in $I$; $q_{i}=\delta \mathcal{G}(\phi)\circ q_{j}$. We must show that there is a unique morphism $q:A\rightarrow G/G'$ such that for every $i\in I$; $q_{i}=\delta p_{i}\circ q$. It is clear that $A$ is an object of Group category and $\{f_{i}:A\rightarrow G_{i}\}_{i\in I}$ is a family of morphisms in $Group$ such that for every $\phi:j\rightarrow i$ in $I$; $f_{i}=\mathcal{G}(\phi)\circ f_{j}$ and for every $i\in I$; $\delta f_{i}=q_{i}$. So there is a unique morphism such as $f:A\rightarrow G$ such that for every $i\in I$; $f_{i}=p_{i}\circ f$. So we consider $q=\delta f$ and because $\delta f_{i}=q_{i}$, therefore the result is obtained.
		Consider the diagram 
		\begin{center}
			\begin{tikzcd}	
				H.Group \ar[r,"\beta"] & Group \\
				I \ar[u,"F"]  \ar[ur,"\beta F"']
			\end{tikzcd}
		\end{center}	
		and suppose that $\underset{I}{lim} \mathcal{F}_{j}$ together with a family $\{t_{i}:\underset{I}{lim} \mathcal{F}_{j}\rightarrow H_{i} \}_{i\in I}$ of morphisms is limit of $I$-diagram $F$. Let $H=\underset{I}{lim} \mathcal{F}_{j}$, then we must show that $\beta(\underset{I}{lim} \mathcal{F}_{j})=H/S$ together with the family $\{\beta t_{i}:H/S\rightarrow H_{i}/S_{i}\}_{i\in I}$ of morphisms is also a limit of $I$-diagram $\beta F$. Where $S=\{h\in H: \beta(h)=e_{\beta(H)}\}$ and $S_{i}=\{h\in H_{i}: \beta(h)=e_{\beta(H_{i})}\}$. \\
		(i)	 Since $L$ is a limit of $I$-diagram $\mathcal{G}$, then for every $\phi:j\rightarrow i$ in $I$; $t_{i}=F(\phi)\circ t_{j}$. So $\beta(t_{i})=\beta F(\phi)\circ\beta(t_{j})$.\\
		(ii) Let $G$ be an object of category $Group$ and $\{s_{i}:G\rightarrow H_{i}/S_{i}\}_{i\in I}$ be a family of morphisms in $Group$ such that for every $\phi:j\rightarrow i$ in $I$; $s_{i}=\beta \mathcal{F}(\phi)\circ s_{j}$. We must show that there is a unique morphism $s:G\rightarrow H/S$ such that for every $i\in I$; $s_{i}=\beta t_{i}\circ s$. is clear that $G$ is an object of $H.Group$ category and $\{g_{i}:G\rightarrow H_{i}\}_{i\in I}$ is a family of morphisms in $H.Group$ such that for every $i\in I$; $\beta(g_{i})=s_{i}$ and for every $\phi:j\rightarrow i$ in $I$; $g_{i}=\mathcal{F}(\phi)\circ g_{j}$. So there is a unique morphism such as $g:G\rightarrow H$ such that for every $i\in I$; $g_{i}=t_{i}\circ g$. So we consider $s=\beta g$ and because $\beta g_{i}=s_{i}$, therefore the result is obtained.\end{proof}}

{\rem	So $\alpha\tiny{\ast}\beta(\underset{j\in I}{lim} \mathcal{F}_{j})
	=\alpha( \underset{j\in I}{lim} \beta(\mathcal{F}_{j}) )
	=\underset{j\in I}{lim}\ \alpha\tiny{\ast}\beta(\mathcal{F}_{j})$ and since $\alpha\tiny{\ast}\beta=\gamma$, then $\gamma$ is continuous and the following diagram is commutative:
	\begin{center}
		\begin{tikzcd}	
			H.Group \ar[r,"\beta"]  \ar[rd]  & Group \ar[d,"\alpha"] \\
			I \ar[dashed,u,"\mathcal{F}"]  \ar[dashed,ur,"\beta\mathcal{F}"'] \ar[dashed,r,"(\alpha\tiny{\ast}\beta)\mathcal{F}=\gamma\mathcal{F}"']	& A.Group
		\end{tikzcd}
	\end{center}}


Let $H$ be a hypergroup, then $E_{H}=\{e\in H|\ x\in e\circ x\cap x\circ e,\ \forall x\in H\}$ and for every $x\in H$; $C_{H}(x)=\{y\in H|\ E_{H}\cap(y\circ x\cap x\circ y)\neq\emptyset\}$. It is clear that $E_{H}\subseteq S_{\gamma_{H}}$ but the equality is not necessarily true. Also $C_{L}(x)=\{y\in H|\ E_{H}\cap y\circ x\neq\emptyset\}$ and $C_{R}(x)=\{y\in H|\ E_{H}\cap x\circ y\neq\emptyset\}$. So $C_{H}(x)=C_{L}(x)\cap C_{R}(x)$ and $y\in C_{L}(x)$ if and only if $x\in C_{R}(x)$. If there is no ambiguity, we write $C(x)$ instead of $C_{H}(x)$.

{\ex
	In any canonical hypergroup $H$ that is not a group, $E_{H}\subseteq S_{\gamma}=S_{\beta}$.}

{\prop
	The hypergroup $H$ is called the regular hypergroup if $E_{H}\neq\emptyset$ and for every $x\in H$; $C_{H}(x)\neq\emptyset$. $H$ is called the strongly regular hypergroup if $E_{H}\neq\emptyset$ and $C_{H}(x)=\{x^{-1}\}$, for every $x\in H$. Clearly if $H$ is strongly regular hypergroup then $|E_{H}|=1$.}

{\cor
	Let $H$ be a strongly regular hypergroup, then $H$ is group if and only if $S_{\beta_{H}}=E_{H}$.}
	{\begin{proof}
		If $|S_{\beta_{H}}=E_{H}|=1$, then $\beta(H)=H/S_{\beta_{H}}\cong H$ and $H$ is group.	\end{proof}	}

{\prop
	Consider the family $\{H_{i}\}_{i\in I}$ of strongly regular hypergroups and suppose that $H_{i}\cap H_{j}=\emptyset$ for every $i,j\in H$. Let $X=\bigcup_{i\in H}H_{i}$, $X\cup\{1\}=\emptyset$. A word on $X$ is any sequence $w=(a_{1},a_{2},...)$ such that $a_{i}\in X\cup \{1\}$ and there is $n\in \mathbb{N}$ such that $a_{i}=1$ for every $i\geq n$. A reduced word on $X$ is a word $w=(a_{1},a_{2},...,a_{n},1,...)$ such that:
	\item[ (i)] For every $i\in I$; $a_{i}\notin E_{H_{i}}$;
	\item[ (ii)] For every $i,j\in I$; if $a_{i}\in H_{j}$, then $a_{i+1}\notin H_{j}$;
	\item[ (iii)] If $a_{k}=1$, then $a_{i}=1$ for every $i\geq k$.\\
	We denote the reduced word $w=(a_{1},a_{2},...,a_{n},1,...)$ by $a_{1}a_{2}...a_{n}$ and $1=(1,1,...)$ is called the empty word.}

\bigskip

Let $\prod^{*}_{i\in I}H_{i}$ be the set of all reduced words on $X=\bigcup_{i\in I}H_{i}$ and $w_{1}=a_{1}a_{2}...a_{n}, w=b_{1}b_{2}...b_{m}\in \prod^{*}_{i\in I}H_{i}$. First of all, we need to distinguish between $a_{1}a_{2}...a_{n}b_{1}b_{2}...b_{m}$ and $a_{1}a_{2}...a_{n}\cdot b_{1}b_{2}...b_{m}$, because the former is a reduced word but the latter is a product of two reduced words that we need to calculate.\\
If $a_{r}\in H_{i_{r}}$($1\leq r\leq n$) and $b_{s}\in H_{j_{s}}$($1\leq s\leq m$), then consider the hyperopration 
$$\cdot:\prod^{*}_{i\in I}H_{i}\times \prod^{*}_{i\in I}H_{i}\rightarrow \mathcal{P}^{*}(\prod^{*}_{i\in I}H_{i})$$ such that $w_{1}\cdot w_{2}=a_{1}a_{2}...a_{n}b_{1}b_{2}...b_{m}$ if $H_{i_{n}}\neq H_{j_{1}}$, and $$w_{1}\cdot w_{2}=\{a_{1}a_{2}...a_{n-1}xb_{2}...b_{m}|\ x\in a_{n}\circ b_{1}\}$$ if $H_{i_{n}}=H_{j_{1}}$, $a_{n}\neq b_{1}^{-1}$. Also $w_{1}\cdot w_{2}=a_{1}a_{2}...a_{n-1}\cdot b_{2}...b_{m}$ if $H_{i_{n}}=H_{j_{1}}$, $a_{n}= b_{1}^{-1}$. Since for every $i\in I$; $H_{i}$ is associative, then $\prod^{*}_{i\in I}H_{i}$ is associative.

{\lem \label{Le}
	Let $\{H_{i}\}_{i\in I}$ be a family of strongly regular hypergroups, then $\prod^{*}_{i\in I}H_{i}$ is strongly regular hypergroup.}
	{\begin{proof}
		Let $w_{1}=a_{1}a_{2}...a_{n}\in\prod^{*}_{i\in I}H_{i}$. Then $w_{1}.\prod^{*}_{i\in I}H_{i}\subseteq\prod^{*}_{i\in I}H_{i}$. If $w_{2}=b_{1}b_{2}...b_{m}\in\prod^{*}_{i\in I}H_{i}$, then consider the reduced word
	\begin{center}	
		$w_{3}=a^{-1}_{n}a^{-1}_{n-1}...a^{-1}_{1}$.
	\end{center}
		If $a^{-1}_{1}$ and $b_{1}$ are not from the same hypergroup, then
	\begin{center}	
		$w_{3}.w_{2}=a^{-1}_{n}a^{-1}_{n-1}...a^{-1}_{1}b_{1}b_{2}...b_{m}$
    \end{center}	
		and $w_{1}.(w_{3}.w_{2})=w_{2}$. If $a^{-1}_{1}$ and $b_{1}$ belong to the same hypergroup and $a^{-1}_{1}\neq b^{-1}_{1}$, then $w_{3}.w_{2}=a_{n}^{-1}a_{n-1}^{-1}...a_{2}^{-1}(a_{1}^{-1}\circ b_{1})b_{2}...b_{m}$ so
		\begin{align*}
		w_{1}.w_{3}.w_{2}&=a_{1}\circ(a_{1}^{-1}\circ b_{1})b_{2}...b_{m}\\
		&=(a_{1}\circ a_{1}^{-1})\circ b_{1}b_{2}...b_{m}\\
		&=b_{1}b_{2}...b_{m}=w_{2}.
	   \end{align*}
		If $a^{-1}_{1}$ and $b_{1}$ belong to the same hypergroup and $a^{-1}_{1}= b^{-1}_{1}$, then $a_{1}=b_{1}$ and $w_{1}.w_{3}.w_{2}=a_{1}a_{2}...a_{n}.a^{-1}_{n}...a_{2}^{-1}a^{-1}_{1}.b_{1}b_{2}...b_{m}=b_{1}b_{2}...b_{m}$. So $w_{1}.\prod^{*}_{i\in I}H_{i}=\prod^{*}_{i\in I}H_{i}$ and $\prod^{*}_{i\in I}H_{i}$ is hypergroup. Since for every $i\in I$, $H_{i}$ is regular, then for every $w\in \prod^{*}_{i\in I}H_{i}$; $|C(w)|=1$ and $(\prod^{*}_{i\in I}H_{i},.)$ is strongly regular hypergroup.\end{proof}}

{\rem
	For every $i\in I$; consider $t_{i}:H_{i}\rightarrow \prod^{*}_{i\in I}H_{i}$ by $t_{i}(x)=x$. Then $t_{i}$ is a good monomorphism of strongly regular hypergroups, because $t_{i}(x\circ y)=\{(z,1,...)|\ z\in x\circ y\}=t_{i}(x)\cdot t_{i}(y)$ and $t_{i}(e_{i})=1$. Also if $w=a_{1}a_{2}...a_{n}$, consider $supp(w)=\{a_{1},a_{2},...,a_{n}\}$.}

{\cor
	If $\{H_{i}\}_{i\in I}$ is a family of polygroups, then $\prod^{*}_{i\in I}H_{i}$ is polygroup.}
	{\begin{proof}
		Since every polygroup is a strongly regular hypergroup, then by Lemma \ref{Le}, $\prod^{*}_{i\in I}H_{i}$ is regular hypergroup, and $E=\{1\}$. Assume that $w_{1}=a_{1}a_{2}...a_{n}$, $w_{2}=b_{1}b_{2}...b_{m}$ and $w_{3}=c_{1}c_{2}...c_{k}$. are reduced words belonging to $\prod^{*}_{i\in I}H_{i}$, such that $w_{1}\in w_{2}.w_{3}$. So $a_{1}a_{2}...a_{n}\in b_{1}b_{2}...b_{m}.c_{1}c_{2}...c_{k}$ and three states will occur.\\
		(i) If $c_{1}$ and $b_{m}$, are not from the same hypergroup, then $w_{1}=w_{2}.w_{3}$ and $w_{1}.w_{3}^{-1}=w_{2}.w_{3}.w_{3}^{-1}$ so $w_{1}.w_{3}^{-1}=w_{2}$. Also $w_{2}^{-1}.w_{1}=w_{2}^{-1}.w_{2}.w_{3}$ so $w_{2}^{-1}.w_{1}=w_{3}$.\\
		(ii) If $c_{1}$ and $b_{m}$, belong to the same hypergroup and $c_{1}=b_{m}^{-1}$, then $c_{1}^{-1}=b_{m}$ and $w_{1}.w_{3}^{-1}\subseteq b_{1}b_{2}...b_{m-1}.c_{1}^{-1}=b_{1}b_{2}...b_{m-1}b_{m}=w_{2}$ so $w_{2}\in w_{1}.w_{3}^{-1}$. Also $w_{2}^{-1}.w_{1}\subseteq b_{m}^{-1}.c_{2}...c_{k}=c_{1}c_{2}...c_{k}=w_{3}$ so $w_{3}\in w_{2}^{-1}.w_{1}$.\\
		(iii) If $c_{1}$ and $b_{m}$, belong to the same hypergroup and $c_{1}\neq b_{m}^{-1}$, then $w_{1}=a_{1}a_{2}...a_{n}\in b_{1}b_{2}...b_{m-1}(b_{m}\circ c_{1})c_{2}...c_{k}=\{b_{1}b_{2}...b_{m-1}xc_{2}...c_{k}; x\in b_{m}\circ c_{1}\}=w_{2}.w_{3}$. Hence, for a member $y\in b_{m}\circ c_{1}$ we have $w_{1}=b_{1}b_{2}...b_{m-1}yc_{2}...c_{k}$.\\
		So $w_{1}.w_{3}^{-1}=b_{1}b_{2}...b_{m-1}yc_{2}...c_{k}.c_{k}^{-1}...c_{2}^{-1}c_{1}^{-1}=b_{1}b_{2}...b_{m-1}y.c_{1}^{-1}=b_{1}b_{2}...b_{m-1}(y\circ c_{1}^{-1})$. Since $H_{i}$'s are polygroup, then $y\in b_{m}\circ c_{1}$ results that $b_{m}\in y\circ c_{1}^{-1}$ and $b_{1}b_{2}...b_{m-1}b_{m}\in b_{1}b_{2}...b_{m-1}(y\circ c_{1}^{-1})$ which means $w_{2}\in w_{1}.w_{3}^{-1}$. Also
		\begin{align*}
		w_{2}^{-1}.w_{1}^{-1}&=b_{m}^{-1}b_{m-1}^{-1}...b_{2}^{-1}b_{1}^{-1}.b_{1}b_{2}...b_{m-1}yc_{2}...c_{k}\\
		&=b_{m}^{-1}.yc_{2}...c_{k}\\
		&=(b_{m}^{-1}\circ y)c_{2}...c_{k}.
		\end{align*}
		Since $y\in b_{m}\circ c_{1}$, then $c_{1}\in b_{m}^{-1}\circ y$ and $c_{1}c_{2}...c_{k}\in (b_{m}^{-1}\circ y)c_{2}...c_{k}$. Hence $w_{3}\in w_{2}^{-1}.w_{1}$ and $\prod^{*}_{i\in I}H_{i}$ is a polygroup.\end{proof}}

{\theo \label{Co}
	$\gamma(\prod^{*}_{i\in I}H_{i})=\sum_{i\in I}\gamma(H_{i})$.}
	{\begin{proof}	Let $H=\prod^{*}_{i\in I}H_{i}$. Since $(H,.)$	 is hypergroup and for every reduced word $w=a_{1}a_{2}...a_{n}$ of $H$, we have $w=a_{1}\textbf{.}a_{2}\textbf{.}...\textbf{.}a_{n}$, then:
		$$\beta(w)=\beta(a_{1})\textbf{.}\beta(a_{2})\textbf{.}...\textbf{.}\beta(a_{n})=\beta(a_{1})\beta(a_{2})...\beta(a_{n})\in \prod^{*}_{i\in I}\beta(H_{i}).$$		 
		If $a\in H\cap H_{i}$, then $\beta(a)=\beta_{i}(a)$. Alao $S_{\beta_{H}}=e_{\beta_{H}}=\beta(w)$, for every $w\in S_{\beta(H)}$ and, we have:
		\begin{align*}
			S_{\beta_{H}} 
			&= \{w\in H; \beta(w).\beta(v)=\beta(v).\beta(w)=\beta(v), \forall v\in H\}\\
			&= \{w\in H; \beta(w.v)=\beta(v.w)=\beta(v), \forall v\in H\} \\
			&= \{a_{1}a_{2}...a_{n}\in H; \beta(a_{1})\beta(a_{2})...\beta(a_{n})=e_{\beta(H)}\}\\
			&= \{a_{1}a_{2}...a_{n}\in H; \beta_{i_{1}}(a_{1})\beta_{i_{2}}(a_{2})...\beta_{i_{n}}(a_{n})=1\}\\
			&= \{a_{1}a_{2}...a_{n}\in H; a_{j}\in S_{i_{j}}, 1\leq j\leq n\}\\
			&= \prod^{*}_{i\in I}S_{i}.	
		\end{align*}
		Since for every $1\leq j\leq n$; $\beta(H_{i_{j}})$ is group and if $a_{j}\in S_{i_{j}}$, then $\beta_{i_{j}}(a_{j})=e_{\beta(H_{i_{j}})}$ and $\beta_{i_{j}}(a_{j})$ can not appear in a reduced word.\\
		Now consider the family of canonical epimorphisms $\{\pi_{i}:H_{i}\rightarrow H_{i}/S_{i}\}_{i\in I}$ and the following map:
		\begin{equation}
		\underset{a_{1}a_{2}...a_{n} \ \longmapsto \ \pi_{i_{1}}(a_{1})\pi_{i_{2}}(a_{2})...\pi_{i_{n}}(a_{n}) \ \ \ (i_{j}\in I)}{\varphi:H\rightarrow \prod^{*}_{i\in I}H_{i}/S_{i}}
		\end{equation}
		Let $w_{1}=a_{1}a_{2}...a_{n}, w_{2}=b_{1}b_{2}...b_{m}\in H$. If $a_{n},b_{1}\notin H_{i}$, for every $i\in I$, then it is clear that $\varphi(w_{1}.w_{2})=\varphi(w_{1}).\varphi(w_{2})$. If there is $k\in I$ such that $a_{n},b_{1}\in H_{k}$ and $a_{n}\neq b_{1}^{-1}$, then $w_{1}.w_{2}=a_{1}a_{2}...a_{n-1}(a_{n}\circ b_{1})b_{2}...b_{m}=\{a_{1}a_{2}...a_{n-1}xb_{2}...b_{m}; x\in a_{n}\circ b_{1}\}\subseteq H$.\\
		Hence $\varphi(w_{1}.w_{2})=\pi_{i_{1}}(a_{1})\pi_{i_{2}}(a_{2})...\pi_{i_{n-1}}(a_{n-1})\pi_{k}(x)\pi_{j_{2}}(b_{2})...\pi_{j_{m}}(b_{m})$ for every $x\in a_{n}\circ b_{1}$.
		But $\pi_{k}(x)=\pi_{k}(a_{n}\circ b_{1})=\pi_{k}(a_{n}).\pi_{k}(b_{1})$, so $\varphi(w_{1}.w_{2})=\varphi(w_{1}).\varphi(w_{2})$.\\
		If $a_{n},b_{1}\in H_{k}$ and $a_{n}=b_{1}^{-1}$, then $w_{1}.w_{2}=a_{1}a_{2}...a_{n-1}.b_{2}...b_{m}$ and $\varphi(w_{1}.w_{2})=\pi_{i_{1}}(a_{1})...\pi_{i_{n-1}}(a_{n-1}).\pi_{j_{2}}(b_{2})...\pi_{j_{m}}(b_{m})$. Since $\beta\subseteq\gamma$, then $\gamma(a_{n}\circ b_{1})=e$ implies that $\beta(a_{n}\circ b_{1})=e$.\\ Hence $\pi_{k}(a_{n}\circ b_{1})=\pi_{k}(a_{n}).\pi_{k}(b_{1})=S_{k}=e_{\beta(H_{k})}$. Therefore
		\begin{align*}
		 \varphi(w_{1}.w_{2})&=\pi_{i_{1}}(a_{1})...\pi_{i_{n-1}}(a_{n-1})\pi_{k}(a_{n}).\pi_{k}(b_{1})\pi_{j_{2}}(b_{2})...\pi_{j_{m}}(b_{m})\\
		 &=\varphi(w_{1}).\varphi(w_{2}).
     	\end{align*} 
		  Also
		\begin{align*}
			ker\varphi&=\{w\in H; \varphi(w)=1\}\\
			&=\{a_{1}a_{2}...a_{n}\in H; \pi_{i_{j}}(a_{j})=S_{i_{j}}\}\\
			&=\{a_{1}a_{2}...a_{n}\in H; \forall i\in[n], \exists j\in I \ni a_{i}\in S_{i_{j}}\}\\
			&=\prod^{*}_{i\in I}S_{i}.
		\end{align*}
		So $H/\prod^{*}_{i\in I}S_{i}\cong\prod^{*}_{i\in I}H_{i}/S_{i}$ and $\beta(\prod^{*}_{i\in I}H_{i})=\prod^{*}_{i\in I}\beta(H_{i})$. Consider the family $\{t_{i}:G_{i}\rightarrow\sum_{i\in I}G_{i}/G_{i}'\}_{i\in I}$ of monomorphisms that $t_{i}(g)=\{x_{i}\}_{i\in I}$ where $x_{i}=g+G_{i}'$ and $x_{j}=G_{j}'$ for $j\neq i$. Let
		\begin{equation}
		\underset{a_{1}a_{2}...a_{n} \ \longmapsto \ \sum_{i=1}^{n}t(a_{i})}{\psi:\prod^{*}_{i\in I}G_{i}\rightarrow \sum_{i\in I}G_{i}/G_{i}'}
		\end{equation}
		If $w_{1}=a_{1}a_{2}...a_{n}, w_{2}=b_{1}b_{2}...b_{m}\in\prod^{*}_{i\in I}G_{i}$, then $w_{1}w_{2}=c_{1}c_{2}...c_{r}$, $r\leq m+n$ and $\psi(w_{1}.w_{2})=\sum_{i=1}^{r}t(c_{i})=\sum_{i=1}^{n}t(a_{i})+\sum_{i=1}^{m}t(b_{i})=\psi(w_{1})+\psi(w_{2})$.\\
		For example, if $\psi:X\ast Y\ast Z\ast W\rightarrow X/X'\oplus Y/Y'\oplus Z/Z'\oplus W/W'$ and $w_{1}=x_{1}y_{1}, w_{2}=z_{1}w_{1}, w_{3}=y_{2}z_{2}x_{2}$ and $w_{4}=y_{1}^{-1}x_{2}w_{2}$ where $x_{i}\in X, y_{i}\in Y, z_{i}\in Z$ and $w_{i}\in W$ for every $i\in\{1,2\}$, then:
		\begin{align*}
			&\psi(w_{1}.w_{2})=\psi(x_{1}y_{1}z_{1}w_{1})= (x_{1}+X',y_{1}+Y',z_{1}+Z',w_{1}+W')\\
			&\psi(w_{1}.w_{3})=\psi(x_{1}y_{1}y_{2}z_{2}x_{2})= \psi(x_{1}yz_{2}x_{2})=(x_{1}x_{2}+X',y+Y',z_{2}+Z',W')\\
			&\psi(w_{1}.w_{4})=\psi(x_{1}x_{2}w_{2})=\psi(xw_{2})= (x+X',Y',Z',w_{2}+W')\\
			&\psi(w_{1}.w_{2}.w_{1}^{-1}.w_{2}^{-1})=\psi(x_{1}y_{1}z_{1}w_{1}y_{1}^{-1}x_{1}^{-1}w_{1}^{-1}z_{1}^{-1})= (X',Y',Z',W').
		\end{align*}
		If $N=\{g_{i}+G_{i}'\}_{i\in I}\in\sum_{i\in I}G_{i}/G_{i}'$, then consider $A=\{g_{i}; \ g_{i}+G_{i}'\in N, g_{i}\notin G_{i}'\}$ and let $w=g_{i_{1}}g_{i_{2}}...g_{i_{n}}$ such that $Supp(w)=A$. Then $\varphi(w)=N$ and $\psi$ is epimorphism. $ker\psi=\{w\in H; x\in Supp(w)\Rightarrow x^{-1}\in Supp(w)\}=\{w\in H; w\in(\prod^{*}_{i\in I}G_{i})'\}=\prod^{*}_{i\in I}G_{i}'$. Let $\prod^{*}_{i\in I}G_{i}=G$, then $\delta(G)=G/G'\cong\sum_{i\in I}G_{i}/G_{i}'=\sum_{i\in I}\delta(G_{i})$. Therefore $\gamma(H)=\delta\tiny{\ast}\beta(H)=\delta(\prod^{*}_{i\in I}\beta(H_{i}))=\sum_{i\in I}\delta\tiny{\ast}\beta(H_{i})=\sum_{i\in I}\gamma(H_{i})$.\end{proof}}

\bigskip

For every $i\in I$; $\bar{t_{i}}:\beta(H_{i})\rightarrow\prod^{*}_{i\in I}\beta(H_{i})$ by $\bar{t_{i}}(\beta(x))=\beta(x)$ is a monomorphism of groups.

{\prop
	A coproduct for the family $\{H_{i}\}_{i\in I}$ of strongly regular hypergroups, is a strongly regular hypergroup $H$ such that $\beta(H)$ is coproduct of the family $\{\beta(H_{i})\}_{i\in I}$ of groups in the group category.}

\bigskip

Equivalently, if $(B,\circ)$ is a regular hypergroup and $\{\psi_{i}:H_{i}\rightarrow B\}_{i\in I}$ is a family of good homomorphisms, then there is a unique homomorphism $\overline{\psi}$ such that the following diagram is commutative:
\begin{center}
	\begin{tikzcd}	
		\beta(H)  \ar[d,"\overline{\psi}"'] & \beta(H_{i}) \ar[dl,"\overline{\psi_{i}}"] \ar[l,"\overline{t_{i}}"']\\
		\beta(B) 
	\end{tikzcd}
\end{center}		
Let $H=\prod^{*}_{i\in I}H_{i}$ and $w=a_{1}a_{2}...a_{n}\in H$, then it is enough to consider $\overline{\psi}:\beta(H)\rightarrow\beta(B)$ by $\overline{\psi}(w)=\beta(\psi_{i_{1}}(a_{1})\circ...\circ\psi_{i_{n}}(a_{n}))=
\overline{\psi}_{i_{1}}(a_{1})...\overline{\psi}_{i_{n}}(a_{n})$.

{\theo
	$\prod^{*}_{i\in I}H_{i}$ is coproduct of the family $\{H_{i}\}_{i\in I}$ of hypergroups.}
	{\begin{proof}
		It can be easily proved, because $\beta(\prod^{*}_{i\in I}H_{i})=\prod^{*}_{i\in I}\beta(H_{i})$, by Theorem \ref{Co}.
	\end{proof}}

{\rem
	If $\{H_{i}\}_{i\in I}$ is a family of abelian strongly regular hypergroups (especially, canoniacl hypergroups), then $\beta(H_{i})=\gamma(H_{i})$, for every $i\in I$. So
	\begin{center}
		$\prod^{*}_{i\in I}\beta(H_{i})=\gamma(\prod^{*}_{i\in I}H_{i})=\sum_{i\in I}\gamma(H_{i})=\sum_{i\in I}\beta(H_{i})=\beta(\sum_{i\in I}H_{i}).$
	\end{center}
	Hence $\sum_{i\in I}H_{i}$  is coproduct of the family $\{H_{i}\}_{i\in I}$ in the corresponding category.}

\bigskip

Since $\delta(\prod^{*}_{i\in I}G_{i})=\sum_{i\in I}\delta(G_{i})$, then $G$ is coproduct of the family $\{G_{i}\}_{i\in I}$ in the category of groups, if $\delta(G)$ is coproduct of the family $\{\delta(G_{i})\}_{i\in I}$ in the category of abelian groups. Therefore we can say that $H$ is coproduct of the family $\{H_{i}\}_{i\in I}$ of strongly regular hypergroups if $\gamma(H)$ is coproduct of the family $\{\gamma(H_{i})\}_{i\in I}$ in the category of abelian groups.

{\theo \label{Co}
	The functors $\delta$ and $\beta$ are cocontinuous.}
	{\begin{proof}
		(a)	 Consider the $I$-diagram $\mathcal{G}:I\rightarrow Group$ and suppose that $G=\underset{j\in I}{colim}\ \mathcal{G}_{i}, \{s_{i}: \mathcal{G}_{i}\rightarrow G\}_{i\in I}$ is colimit of $I$-diagram $\mathcal{G}$.\\
		(i) Since for every $\varphi:i\rightarrow j$ in $I$; $s_{i}=s_{j}\circ\mathcal{G}(\varphi)$, then $\delta(s_{i})=\delta(s_{j})\circ\delta\mathcal{G}(\varphi)$. (ii) Let $Z$ be an abelian group and $\{r_{i}; \delta\mathcal{G}_{i}\rightarrow Z\}_{i\in I}$ be a family of homomorphisms in A.Group, such that for every $\varphi:j\rightarrow i$ in $I$; $r_{i}=r_{j}\circ\delta\mathcal{G}(\varphi)$. Hence there is a family $\{r_{i}'; \mathcal{G}_{i}\rightarrow Z\}_{i\in I}$ in Group, such that $\alpha(r_{i}')=r_{i}$ for every $i\in I$. So consider $r_{i}'=r_{j}'\circ\mathcal{G}(\varphi)$ for every $i\in I$. Since $G$ is colimit of $I$-diagram $\mathcal{G}$, then there is a unique morphism $t:G\rightarrow Z$  such that $r_{i}'=t\circ s_{i}$. Therefore, there is a unique morphism $\delta t:\delta(G)\rightarrow \delta(Z)$ such that $\delta(r_{i}')=r_{i}=\delta(t)\circ\delta(s_{i})$. Hence $\delta(\underset{j\in I}{colim}\ \mathcal{G}_{i})=\underset{j\in I}{colim}\ \delta\mathcal{G}_{i}$.\\
		(b) Consider the $I$-diagram $\mathcal{F}:I\rightarrow H.Group$ and suppose that
		\begin{center}
			$H=\underset{j\in I}{colim}\ \mathcal{F}_{i} \ , \{s_{i}: \mathcal{F}_{i}\rightarrow H\}_{i\in I}$
		\end{center}
		is colimit of $I$-diagram $\mathcal{F}$.
		(i) Since for every $\varphi:i\rightarrow j$ in $I$; $s_{i}=s_{j}\circ\mathcal{F}(\varphi)$, then $\beta(s_{i})=\beta(s_{j})\circ\beta\mathcal{F}(\varphi)$.(ii) Let $G$ be a group and $\{r_{i}; \beta\mathcal{F}_{i}\rightarrow G\}_{i\in I}$ be a family of homomorphisms in Group, such that for every $\varphi:i\rightarrow j$ in $I$, $r_{i}=r_{j}\circ\beta\mathcal{F}(\varphi)$. Hence there is a family $\{r_{i}'; \mathcal{F}_{i}\rightarrow G\}_{i\in I}$ in H.Group, such that $\beta(r_{i}')=r_{i}$ for every $i\in I$. So consider $r_{i}'=r_{j}'\circ\mathcal{F}(\varphi)$ for every $i\in I$. Since $H$ is colimit of $I$-diagram $\mathcal{F}$, then there is a unique morphism $t:H\rightarrow G$  such that $r_{i}'=t\circ s_{i}$. Therefore, there is a unique morphism $\beta t:\beta(H)\rightarrow G$ such that $\beta(r_{i}')=r_{i}=\beta(t)\circ\beta(s_{i})$. Hence $\beta(\underset{j\in I}{colim}\ \mathcal{F}_{i})=\underset{j\in I}{colim}\ \delta\mathcal{F}_{i}$.\end{proof}}

{\rem
	Therefore $\delta\tiny{\ast}\beta(\underset{j\in I}{colim}\  \mathcal{F}_{j})=\delta(\underset{j\in I}{colim}\ \beta(\mathcal{F}_{j}))=\underset{j\in I}{colim}\  \delta\tiny{\ast}\beta(\mathcal{F}_{j})$ and since $\gamma=\delta\tiny{\ast}\beta$, then $\gamma$ is cocontinuous.}

\section{ Conclusions}
It is proved that for a closed subhypergroup $K$ of the hypergroup $H$, which contains $S_{\beta}$ (resp. $S_{\gamma}$), the congruence relation modulo $K$ is a strongly regular relation and $H/K$ is a group (resp. abelian group).\\
Additionally, we define the reduced words on strongly regular hypergroups and examined the continuity and cocontinuity of the functor $\gamma^{\ast}$. Also, the behavior of free products under fundamental relations $\beta^{\ast}$ and $\gamma^{\ast}$ was investigated. After that, the relationship between $S(H)$ and $S(H/K) $ was established when $K$ is a canonical subhypergroup of canonical hypergroup $H$. Since every canonical hypergroup is a strongly regular hypergroup, we can now combine the results of these two sections.\\
To this end, assume that $\{H_{i}\}_{i\in I}$ and $\{K_{i}\}_{i\in I}$ are families of canonical hypergroups, where each $K_{i}$ is a canonical subhypergroup of $H_{i}$. In this case, $\{H_{i}/K_{i}\}_{i\in I}$ is a family of canonical hypergroups on which reduced words and coproduct $\prod^{*}_{i\in I}H_{i}/K_{i}$ can be defined and the following equations can be claimed:
\begin{center} 
	$\prod^{*}_{i\in I}H_{i}/K_{i}=\dfrac{\prod^{*}_{i\in I}H_{i}}{\prod^{*}_{i\in I}K_{i}}$,
	$\boldsymbol{\beta}\big{(}\prod^{*}_{i\in I}H_{i}/K_{i}\big{)}=\dfrac{\prod^{*}_{i\in I}H_{i}}{\prod^{*}_{i\in I}S_{\beta_{i}}\circ \prod^{*}_{i\in I}K_{i}}$,
	$\boldsymbol{\gamma}\big{(}\prod^{*}_{i\in I}H_{i}/K_{i}\big{)}=\prod^{*}_{i\in I}\dfrac{H_{i}}{H'_{i}\circ K_{i}}$.
\end{center}
Now the question arises whether the closure condition can be removed? Also, is there a one-to-one correspondence between the set $SR(H)$ and the set of all normal subhypergroups of $H$, which contains $S_{\beta}$? And can the formula $\boldsymbol{\rho}(H/K)=\dfrac{H}{s_{\rho}\circ K}$ be extended to arbitrary hypergroup $H$, that is not necessarily canonical?\\
Also, we hope that this paper provide a suitable tools to deep study on (strongly) regular relations on hyperrings and hypermodules.

\bigskip\bigskip\bigskip
\renewcommand{\bibsection}{\begin{flushright}\Large
{
REFERENCES}\\
\rule{8cm}{0.4pt}\\[0.8cm]
\end{flushright}}

\bigskip\bigskip


\begin{flushleft}
\rule{8cm}{0.4pt}\\
\end{flushleft}

{
\sloppy
\noindent
Behnam Afshar

\noindent 
Department of Mathematics, Statistics and Computer Science

\noindent 
College of Science

\noindent 
University of Tehran

\noindent 
Tehran  (Iran)

\noindent
e-mail: behnamafshar@ut.ac.ir \ \ \ ORCID: 0009-0006-7095-5795
}

\bigskip
\bigskip

{
\sloppy
\noindent
Reza Ameri\ \ \ \ \ (corresponding author)

\noindent 
Department of Mathematics, Statistics and Computer Science

\noindent 
College of Science

\noindent 
University of Tehran

\noindent 
Tehran  (Iran)

\noindent
e-mail: rameri@ut.ac.ir \ \ \ ORCID: 0000-0001-5760-1788

}

\hypertarget{last_page}{}\label{last_page}

\end{document}